# A Singular Miquel Configuration and the Miquel Direct Similarity

**Christopher J Bradley**

**Abstract:** Let ABC be a triangle with P on AB, and let circle APC meet BC at Q and circle BPC meet CA at R, then the special Miquel configuration is when P, Q, R are the operative points on the sides. We show that in this case if S is the Miquel point then ASQ and BSR are straight lines. In the last part we investigate the direct similarity between ABC and DEF the triangle formed by the centre of the Miquel circles.

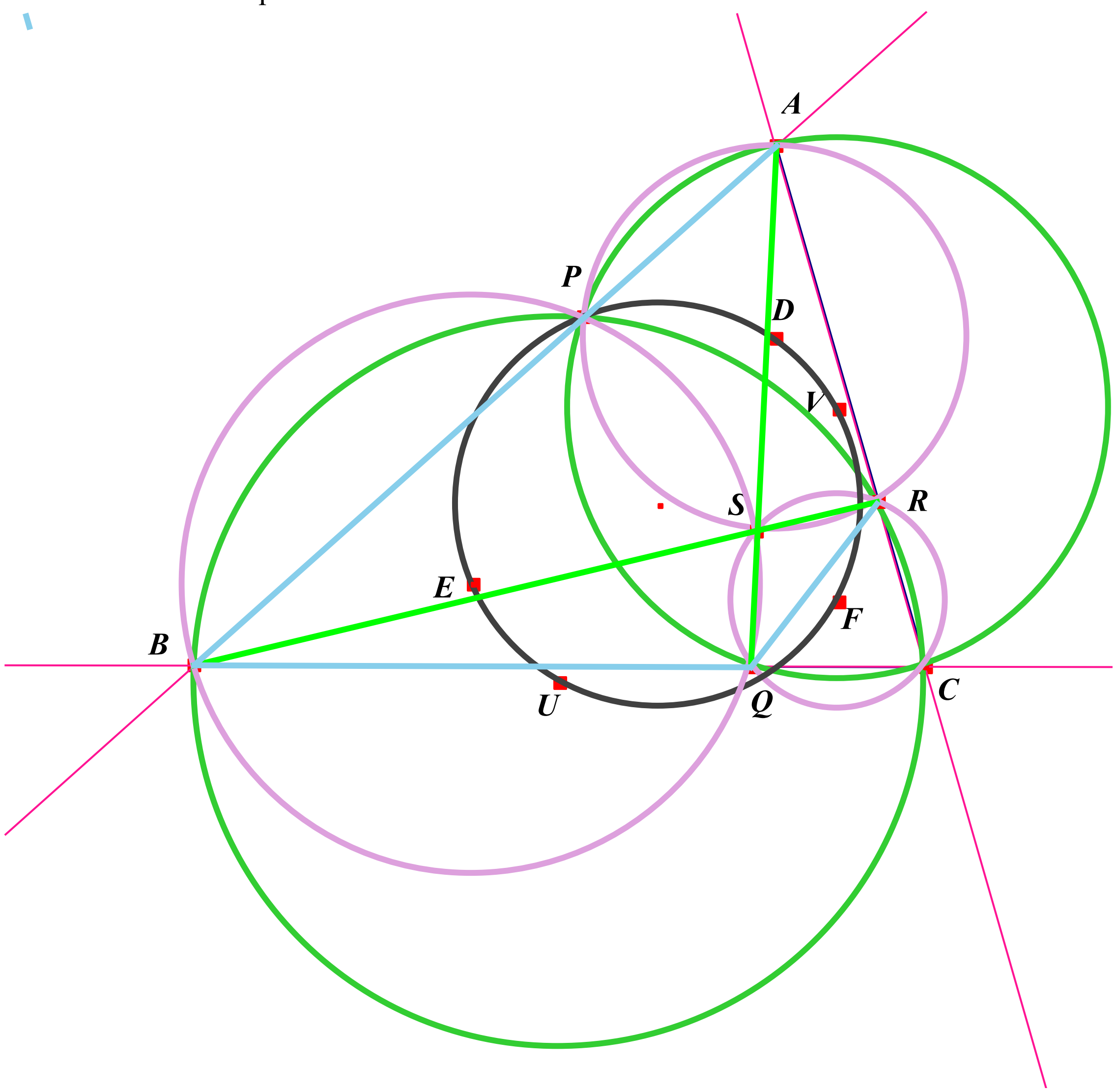

**Fig. 1**
**The singular Miquel configuration**

## 1. Introduction

A Miquel configuration is one in which points P, Q, R are chosen on the sides AB, BC, CA respectively (but not at the vertices A, B, C) and the circles ARP, BPQ and CQR are drawn, illustrating that these three circles always pass through a common point S. The configuration we have called a *singular Miquel configuration* is that once P has been chosen arbitrarily (but not at B or C), then Q and R are chosen so that Q lies on circle CAP and R lies on the circle BCP. We prove in this paper that in these circumstances ASQ and BSR are straight lines and that the centres U, V, D, E, F of the five circles BCP, CAP, ARP, BPQ, CQR respectively do themselves lie on a circle passing through P. Fig. 1 shows a case in which P lies between A and B, but the results still holds when P is external and the algebra we present is independent of the position of P on AB. We also show that triangle DEF is directly similar to ABC and as this is true in any Miquel configuration we give a proof when P, Q and R are arbitrary. In the analysis of the singular Miquel configuration we use areal co-ordinates in which ABC is the reference triangle. For the study of the direct similarity we use Cartesian co-ordinates.

## 2. Finding Q and R from P

Suppose P is chosen to have co-ordinates (n, (1 – n), 0), $n \neq 0, 1$. The general equation of a circle in areal co-ordinates is

$$a^2yz + b^2zx + c^2xy - (x+y+z)(ux+vy+wz) = 0. \qquad (2.1)$$

Here u, v, w are determined by substituting the co-ordinates of three points lying on the circle. Using this method repeatedly we find the equations of CAP and BCP are respectively are

$$c^2ny^2 + (c^2n - a^2)yz - b^2zx + c^2(n-1)xy = 0 \qquad (2.2)$$

and

$$c^2(n-1)x^2 + a^2yz + \left(b^2 + c^2(n-1)\right)zx + c^2nxy = 0. \qquad (2.3)$$

The co-ordinates of Q are therefore $(0, a^2 - c^2n, c^2n)$ and the co-ordinates of R are $(b^2 + c^2(n-1), 0, c^2(1-n))$. The point S we define to be AQ^BR and has co-ordinates $(n(b^2 + c^2(n-1)),$ $(1-n)(a^2 - c^2n), c^2n(1-n))$.

## 3. The three Miquel circles

The circles ARP. BPQ have equations

$$c^2ny^2 + \left(b^2 + c^2(n-1)\right)z^2 - \left(a^2 - b^2 + c^2(1-2n)\right)yz + c^2(n-1)zx + c^2(n-1)xy = 0, \qquad (3.1)$$

$$c^2(n-1)x^2 + (c^2n - a^2)z^2 + c^2nyz - \left(a^2 - b^2 + c^2(1-2n)\right)zx + c^2nxy = 0. \qquad (3.2)$$

and it may be checked that both these circles pass through S.

Finally, the circle CQR has equation

$$c^2(n-1)x^2 - c^2ny^2 - (c^2n - a^2)yz + \left(b^2 + c^2(n-1)\right)zx = 0, \tag{3.3}$$

and, of course, it too passes through S.

### 4. The centres of the 5 circles

If a circle is expressed in the form of equation (2.1), then, see Bradley [1], its centre has co-ordinates x, y, z), where

$$\begin{aligned} x &= -\left(\tfrac{1}{4}\right)(a^4 - a^2(b^2 + c^2 - 2u + v + w) + (b^2 - c^2)(w - v)), \\ y &= \left(\tfrac{1}{4}\right)(a^2(b^2 + u - w) - b^4 + b^2(c^2 + u - 2v + w) + c^2(w - u)), \\ z &= \left(\tfrac{1}{4}\right)(a^2(c^2 + u - v) + b^2(c^2 - u + v) - c^4 + c^2(u + v - 2w)). \end{aligned} \tag{4.1}$$

Using Equation (4.1) and values of u, v, w from the working in Section 2 we find the co-ordinates of the centre of circle CAP are

$$\begin{aligned} x &= -\left(\tfrac{1}{4}\right)(a^4 - a^2\left(b^2 + c^2(n+1)\right) - c^2n(b^2 - c^2)), \\ y &= \left(\tfrac{1}{4}\right)(b^2\left(a^2 - b^2 + c^2(1-2n)\right)), \\ z &= -\left(\tfrac{1}{4}\right)(c^2\left(a^2(n-1) - b^2(n+1) + c^2(1-n)\right)). \end{aligned} \tag{4.2}$$

and the co-ordinates of the centre of circle BCP are

$$\begin{aligned} x &= -\left(\tfrac{1}{4}\right)(a^2\left(a^2 - b^2 + c^2(1-2n)\right)), \\ y &= \left(\tfrac{1}{4}\right)(a^2\left(b^2 + c^2(1-n)\right) - b^4 + b^2c^2(2-n) + c^4(n-1)), \\ z &= \left(\tfrac{1}{4}\right)(c^2\left(n(b^2 - c^2) - a^2(n-2)\right)). \end{aligned} \tag{4.3}$$

The co-ordinates of the centre of circle ARP are

$$\begin{aligned} x &= -\left(\tfrac{1}{4}\right)(a^4 - 2a^2(b^2 + c^2n) + (b^2 - c^2)^2), \\ y &= \left(\tfrac{1}{4}\right)\left(c^2(1-n)(a^2 + b^2 - c^2)\right), \\ z &= \left(\tfrac{1}{4}\right)\left(c^2(1-n)(a^2 - b^2 + c^2)\right). \end{aligned} \tag{4.4}$$

The co-ordinates of the centre of circle BPQ are

$$x = \left(\frac{1}{4}\right)\left(c^2 n(a^2 + b^2 - c^2)\right),$$
$$y = -\left(\frac{1}{4}\right)\left(a^4 - 2a^2(b^2 + c^2) + b^4 + 2b^2c^2(n-1) + c^4\right), \qquad (4.5)$$
$$z = -\left(\frac{1}{4}\right)\left(c^2 n(a^2 - b^2 - c^2)\right).$$

Finally the co-ordinates of the centre of the circle CQR are

$$x = -\left(\frac{1}{4}\right)\left(a^4 - a^2\left(b^2 + c^2(3n-1)\right) - c^2 n(b^2 - c^2)\right),$$
$$y = \left(\frac{1}{4}\right)\left(a^2\left(b^2 + c^2(1-n)\right) - b^4 + b^2c^2(2-3n) + c^4(n-1)\right), \qquad (4.6)$$
$$z = \left(\frac{1}{2}\right)\left(c^2\left(b^2 n - a^2(n-1)\right)\right).$$

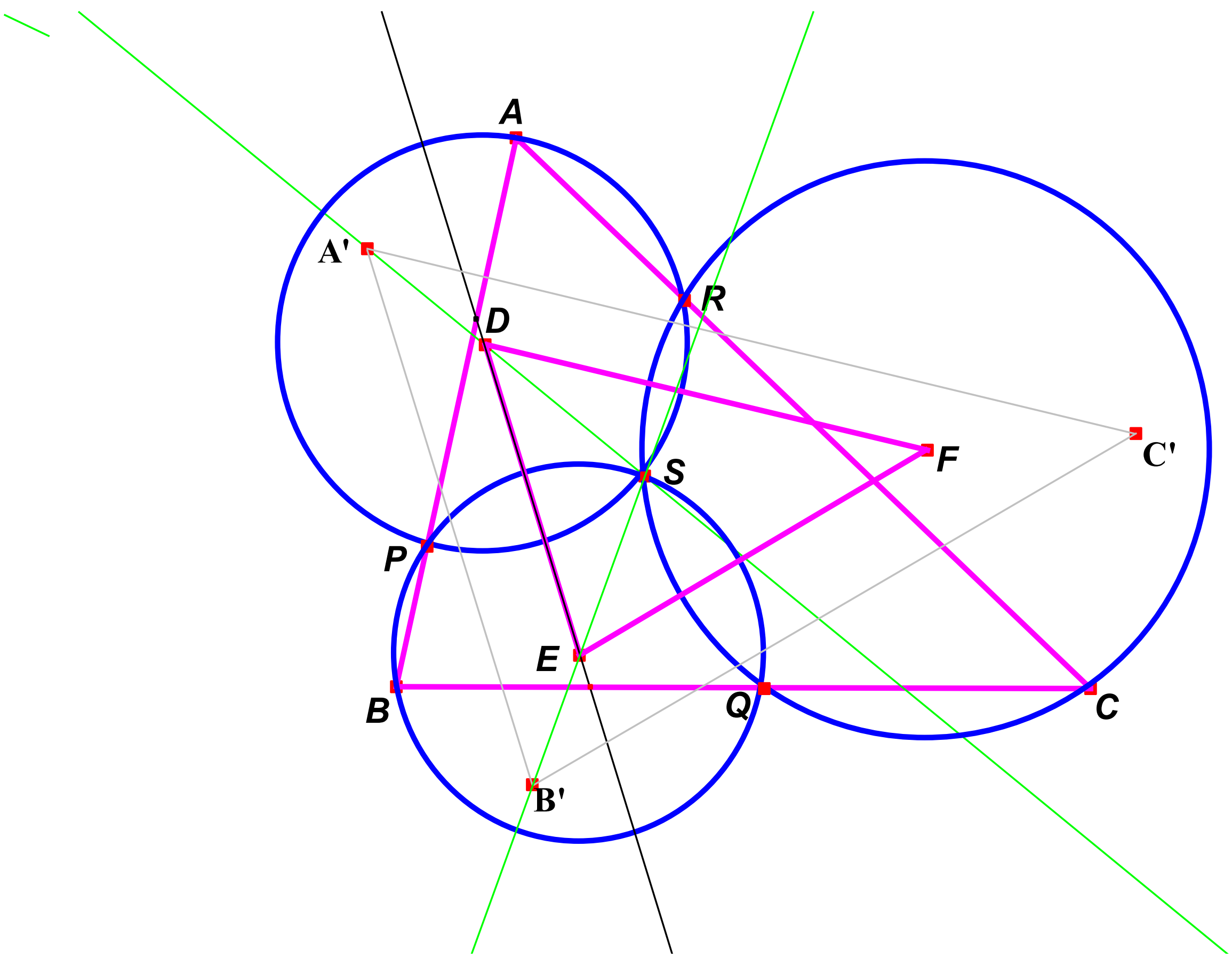

**Fig. 2**
**The direct similarity**

The last three of these points was used to compute the circumcircle of triangle DEF and it was then checked that it does actually contain U, V and P. The circumcircle has a cumbersome equation which is

$$c^2(n-1)\left(a^2\left(b^2+c^2(1-n)\right)-b^4+b^2c^2(2-n)+c^4(n-1)\right)x^2$$
$$+c^2n\left(a^4-a^2\left(b^2+c^2(n+1)\right)-c^2n(b^2-c^2)\right)y^2$$
$$+c^2(a^2+b^2-c^2)\left(a^2(n-1)-n\left(b^2+c^2(n-1)\right)\right)z^2 \quad (4.7)$$
$$-\left(a^6-a^4\left(2b^2+c^2(2n+1)\right)+a^2\left(b^4+b^2c^2(n-1)+c^4n(2n+1)\right)\right.$$
$$\left.+c^2n(b^2-c^2)\left(b^2+c^2(2n-1)\right)\right)yz$$
$$-\left(a^4\left(b^2+c^2(1-n)\right)-a^2\left(2b^4+b^2c^2n+2c^4n(1-n)\right)+b^6+b^4c^2(2n-3)\right.$$
$$\left.+b^2c^4(n-1)(2n-3)+c^6(1-n)(2n-1)\right)zx$$
$$+(c^2\left(a^4(n-1)+a^2\left(b^2-c^2(2n^2-n-1)\right)-n\left(b^4+b^2c^2(2n-3)+2c^4(1-n)\right)\right)xy$$
$$=0.$$

## 5. The direct similarity when P, Q, R are arbitrary

For the more general Miquel configuration, when P, Q, R are chosen as arbitrary points on AB, BC, CA respectively, it is easier to use Cartesian co-ordinates. An economically chosen set of parameters involve co-ordinates (0, 0) for A, (2, 2v) for B and (2, 2w) for C, where without loss of generality we take w > v. Points P, Q, R on AB, BC, CA respectively are now assigned co-ordinates (k, kv), (2, 2u), (h, hw) respectively, where $u \neq v, w$; $h, k \neq 0,2$. See Fig. 2.

The equations of the circles BPQ, CQR, ARP may be obtained using standard methods involving 4 x 4 determinants and are respectively

$$x^2+y^2+(2(uv-1)-k(v^2+1))x-2(u+v)y+2k(v^2+1)=0, \quad (5.1)$$
$$x^2+y^2+(2(uw-1)-h(w^2+1))x-2(u+w)y+2h(w^2+1)= \quad (5.2)$$
$$(v-w)(x^2+y^2)+(kw(v^2+1)-hv(w^2+1))x+(h(w^2+1)-k(v^2+1))y=0. \quad (5.3)$$

It may be checked that these three circles have a common point S. We do not record its co-ordinates as they are complicated and are not needed in what follows. The co-ordinates of their centres can now be determined and are

D (centre of ARP): $\{1/(2(v-w))\}(hv(w^2+1)-kw(v^2+1), k(v^2+1)-h(w^2+1))$,

E (centre of BPQ): $(\tfrac{1}{2}(k(v^2+1)-2(uv-1)), u+v)$,

F (centre of CQR) : $(\tfrac{1}{2}(hw^2+1)-2(uw-1)), u+w)$.

Calculations may now be carried out to show that DE/EF = AB/BC = $\{\sqrt{(v^2+1)}\}/(w-v)$, and DF/DE = AC/AB = $\sqrt{\{(w^2+1)/(v^2+1)\}}$ and this establishes the similarity between triangle *XYZ* and ABC. It may be shown that triangle ABC may be moved by appropriate rotation to A'B'C' and enlargement (reduction) about S as centre of direct similarity to coincide with triangle DEF.

Flat 4,
Terrill Court,
12-14 Apsley Road
BRISTOL BS8 2SP